\documentclass[sts]{imsart}
\usepackage{graphicx} %

\usepackage{amsfonts,amssymb,float}

\DeclareSymbolFont{sfletters}{OML}{cmbrm}{m}{it}

\DeclareMathSymbol{\salpha}{\mathord}{sfletters}{"0B}

\usepackage{mathrsfs}
\RequirePackage{natbib}
\RequirePackage[colorlinks,citecolor=blue,urlcolor=blue]{hyperref}

\usepackage{xcite}
\makeatletter
\newcommand*{\addFileDependency}[1]{%
  \typeout{(#1)}%
  \@addtofilelist{#1}%
  \IfFileExists{#1}{}{\typeout{No file #1.}}%
}
\makeatother

\usepackage{mathtools}
\usepackage{subfigure} 
\RequirePackage{natbib}
\providecommand{\noopsort}[1]{}
\usepackage{comment}
\usepackage{enumerate}
\usepackage{eufrak}

\usepackage{array,varwidth}

\RequirePackage[OT1]{fontenc}
\RequirePackage{amsthm,amsmath}

\startlocaldefs

\numberwithin{equation}{section}
\theoremstyle{plain}
\newtheorem{theorem}{Theorem}[section]
\newtheorem{lemma}{Lemma}[section]
\newtheorem{corollary}{Corollary}[section]

\newtheorem{proposition}{Proposition}[section]
\newtheorem{assumption}{Assumption}[section]

\newcommand{\R}{\mathbb{R}}

\renewcommand{\hat}[1]{\widehat{#1}}

\newcommand{\I}{\textsc{I}}
\newcommand{\II}{\textsc{II}}
\newcommand{\III}{\textsc{III}}

\newcommand{\dK}{d_{\textup{K}}}

\newcommand{\J}{\mathcal{J}}

\renewcommand{\P}{\mathbb{P}}
\newcommand{\E}{\mathbb{E}}

\renewcommand{\L}{\mathcal{L}}

\newcommand{\e}{\epsilon}

\newcommand{\tr}{\operatorname{tr}}

\newcommand{\var}{\operatorname{var}}

\newcommand{\ttop}{^{\top}}

\newcommand{\op}{_{\textup{op}}}

\newcommand{\ts}{\textstyle}

\endlocaldefs

\def\references{\bibliography{lower_tail_bib.bib}}

\begin{document}

\begin{frontmatter}
\title{A Sharp Lower-tail Bound for Gaussian Maxima with Application to Bootstrap Methods in High Dimensions}

\begin{aug}
\vspace{0.2cm}
\author{\fnms{Miles E.} \snm{Lopes}\thanksref{t1}\ead[label=e1]{melopes@ucdavis.edu}}
\and
\author{\fnms{Junwen} \snm{Yao}\ead[label=e2]{jwyao@ucdavis.edu}}

\thankstext{t1}{Supported in part by NSF grant DMS 1613218 and 1915786}

\affiliation{University of California, Davis}

\end{aug}

\begin{abstract}
Although there is an extensive literature on the maxima of Gaussian processes, there are relatively few non-asymptotic bounds on their lower-tail probabilities. 
The aim of this paper is to develop such a bound, while also allowing for many types of dependence. Let
 $(\xi_1,\dots,\xi_N)$ be a centered Gaussian vector with standardized entries, whose correlation matrix $R$ satisfies $\max_{i\neq j} R_{ij}\leq \rho_0$ for some constant $\rho_0\in (0,1)$. Then, for any $\e_0\in(0,\sqrt{1-\rho_0})$, we establish an upper bound on the probability  $\P(\max_{1\leq j\leq N} \xi_j\leq \e_0\sqrt{2\log(N)})$ in terms of $(\rho_0,\e_0,N)$.
The bound is also sharp, in the sense that it is attained up to a constant, independent of $N$. Next, we apply this result in the context of high-dimensional statistics, where we simplify and weaken conditions that have recently been used to establish near-parametric rates of bootstrap approximation. Lastly, an interesting aspect of this application is that it makes use of recent refinements of Bourgain and Tzafriri's ``restricted invertibility principle''.
\end{abstract}



\end{frontmatter}

\section{Introduction} 

The maxima of Gaussian processes play an essential role in many aspects of probability and statistics, and the literature describing them is highly developed~\citep[e.g.][]{Leadbetter:1983,Adler:1990,Lifshits:1995,Ledoux:Talagrand:2013, Talagrand:2014}. Within this area,  a variety of questions are related to showing that the maximum of a process is unlikely to deviate far above, or below, its mean. 
However, in comparison to the set of tools for handling the upper tail of a maximum, there are relatively few approaches for the lower tail. (Additional commentary related to this distinction may be found in~\cite[p.viii]{Talagrand:2014}~\cite[Sec.4.2]{Li:Shao:2001}~\cite[Sec.18]{Lifshits:1995}.)

In this paper, our goal is to derive lower-tail bounds for Gaussian maxima that are motivated by statistical applications involving bootstrap methods in high dimensions. We desire bounds that are general enough to handle many types of correlation structures, yet also precise enough to yield explicit rates of convergence in distributional approximation results. 

To describe the setting of our lower-tail bounds, let $\xi=(\xi_1,\dots,\xi_N)$ be Gaussian vector with a correlation matrix $R\in\R^{N\times N}$, as well as $\E(\xi_j)=0$ and $\var(\xi_j)=1$ for all $1\leq j\leq N$. We will consider the situation where there is a fixed constant $\rho_0\in (0,1)$ such that
\begin{equation}\label{eqn:assumcor}
\max_{i\neq j}R_{ij}\, \leq \, \rho_0.
\end{equation}
In addition, we will consider some relaxations of this condition, where $\rho_0=1$, or where only subset of the off-diagonal entries of $R$ are bounded above by a given constant. (See Corollary~\ref{cor:simple} and Theorem~\ref{thm:key}.)
It is also worth noting that~\eqref{eqn:assumcor} does not require $R$ to be invertible. Letting the maximum entry of $\xi$ be denoted as 
$$M_N(\xi) \, = \, \max_{1\leq j\leq N} \xi_j,$$
we seek non-asymptotic upper bounds on the probability
$\P(M_N(\xi) \leq t)$, where $t$ is a suitable point in the lower tail. %

\subsection{Background }
We now briefly review some leading results on lower-tail bounds for $M_N(\xi)$.
Under the preceding conditions, the well-known concentration inequality for Lipschitz functions of Gaussian vectors implies that for any $s>0$,
\begin{equation}\label{eqn:BTS}
\P\Big( M_N(\xi) \ \leq \ \text{med}(M_N(\xi))-s\Big) \, \leq \, e^{-s^2/2},
\end{equation}
where $\text{med}(\cdot)$ is any median~\citep{Sudakov:1974,Borell:1975}.
Although this bound is broadly applicable, it can fail to describe lower-tail probabilities smaller than $\mathcal{O}(N^{-1})$.  To see this, consider using~\eqref{eqn:BTS} to bound the probability $\P(M_N(\xi)\leq \delta_0\text{med}(M_N(\xi)))$ for some fixed $\delta_0\in(0,1)$. If the entries of $\xi$ are independent, then the standard fact that $\text{med}(M_N(\xi))=\sqrt{2\log(N)}(1+o(1))$ as $N\to\infty$
implies that~\eqref{eqn:BTS} cannot give a bound better than $\mathcal{O}(N^{-1})$ in this case. 
Furthermore, such a bound is far too large. In fact, when the entries of $\xi$ are independent, this probability is $\mathcal{O}(\exp(-a N^b))$, for some fixed constants $a,b>0$ that may depend on $\delta_0$ but not $N$~\citep[cf.][Section 2]{Schechtman}.
More recently, an important result of Paouris and Valettas~\cite[]{Paouris:2018} went beyond~\eqref{eqn:BTS}, showing that the inequality
\begin{equation}\label{eqn:paouris}
\P\Big(M_N(\xi)\leq \text{med}(M_N(\xi))-s\Big) \ \leq \ \ts\frac{1}{2}\exp\Big(-\ts\frac{\pi}{1024}\frac{s^2}{\var(M_N(\xi))}\Big)
\end{equation}
holds for any $s>0$, which can improve upon~\eqref{eqn:BTS} in situations where $\var(M_N(\xi))$ is small.
For further variations and related results, we refer to the papers~\cite{Paouris:Euclidean,Tanguy:2019esaim,Valettas:2017}.)\\[-0.2cm]

Despite the progress achieved by the bound~\eqref{eqn:paouris}, it can still be quite challenging to obtain precise control on the variance $\var(M_N(\xi))$ in the exponent.
Indeed, it is often the case that bounds on $\var(M_N(\xi))$ have implicit dependence on the correlation matrix $R$ that is difficult to quantify, and such bounds usually involve constants that are unspecified or conservative.  (Note that since $\var(M_N(\xi))$  appears in the exponent, the constant in a variance bound will typically affect the rate of the tail bound.) %
 We refer to the books~\citep{Massart:2013,Chatterjee:2014} for further background on variance bounds related to $M_N(\xi)$.\\[-0.2cm]

\noindent\textbf{Contributions.} With regard to the considerations above, a few other works have developed lower-tail bounds for $M_N(\xi)$ that provide insight into the role of the correlation structure~\citep{Hartigan:2014,Ding:2015,Tanguy:2015}. However, a limitation shared by all of these works is that they do not explicitly quantify the rate at which probabilities of the form $\P(M_N(\xi)\leq \e_0\sqrt{2\log(N)})$ decrease with respect to $N$ for a fixed value $\e_0\in (0,\sqrt{1-\rho_0})$.
For instance, this limitation can arise from unspecified constants in the exponents of the bounds. See also the discussion after Theorem~\ref{THM:MAIN} below for more details.

In light of this, an important contribution of our lower-tail bounds is that they provide rates with explicit constants in their exponents.
Also, our work goes further by showing that the rates are \emph{sharp} (as in~\eqref{eqn:mainbound} and~\eqref{eqn:mainlower} of Theorem~\ref{THM:MAIN}). Moreover, the proof techniques used here are entirely different from those used in the previously mentioned works. In particular, we extend and apply our lower-tail bounds by leveraging recent refinements of Bourgain and Tzafriri's \emph{restricted invertibility principle} (Section~\ref{sec:ri}), which has been a topic of substantial interest in contemporary mathematics. Conventionally, this principle is understood as a functional-analytic result that guarantees the existence of special submatrices within large matrices. Hence, our use of this result to serve the quite different purpose of enhancing tail bounds may be viewed as another notable aspect of our work.

In addition to these contributions, Theorem~\ref{THM:BOOT} in Section~\ref{sec:app} shows how our lower-tail bounds can be applied in the context of high-dimensional statistics, in order to simplify and weaken conditions that are sufficient for near-parametric rates of bootstrap approximation. Specifically, Theorem~\ref{THM:BOOT} advances the state of the art on such results for ``max statistics'' in settings where the data satisfy a condition known as ``variance decay''. Also, a second application of our results has recently been developed by~\citep{Yi:2021} in connection with stochastic PDEs. Brief descriptions of both applications are given below.

\subsection{Applications}\label{sec:appintro}
\noindent \textbf{Bootstrap methods in high dimensions.} In recent years, inference problems related to max statistics have attracted significant attention in the high-dimensional statistics literature. The prototypical example of a max statistic is the coordinate-wise maximum of a sum of random vectors $X_1,\dots,X_n\in\R^p$, say
\begin{equation}\label{eqn:maxstatdef}
T=\max_{1\leq j\leq p}\frac{1}{\sqrt n}\displaystyle\sum_{i=1}^n X_{ij}.
\end{equation}
Statistics of this type arise often in the construction of hypothesis tests and confidence regions, and consequently, the performance of many inference procedures is determined by how accurately the distribution $\mathcal{L}(T)$ can be approximated. %

Bootstrap methods are a general approach to construct such approximations in a data-driven way, and they are designed to generate a random variable $T^*$ whose conditional distribution given the data, denoted $\mathcal{L}(T^*|X)$, is close to $\mathcal{L}(T)$. In this context, accuracy is commonly measured with respect to the Kolmogorov metric $d_{\textup{K}}$,
 and the challenge of developing non-asymptotic upper bounds on $d_{\textup{K}}(\mathcal{L}(T^*|X),\mathcal{L}(T))$ has stimulated an active line of work~\citep[e.g.][]{Chernozhukov:2019,Chernozhukov:2020,Deng:2020,Lopes:2020,Lopes2:2020}.
As an application of our lower-tail bounds, we will consider a recent bootstrap approximation result of this type from~\citep{Lopes:2020}, and we will improve upon it by showing that it holds under assumptions that are simpler and less restrictive.\\[-0.2cm]

\noindent \textbf{Macroscopic properties of solutions to stochastic PDEs.}
In the study of stochastic partial differential equations, it is of interest to determine whether or not solutions exhibit high peaks over large regions at different scales. Solutions having this property are said to be ``multifractal at macroscopic scales''. In order to demonstrate that a solution has this property in  a precise sense, it is necessary to analyze regions where a solution rises above a certain height (exceedance sets), and quantify the ``macroscopic Hausdorff dimension'' of such regions. 

During the last few years, a growing number of results have demonstrated that solutions to fundamental stochastic PDEs possess this multifractical property~\citep[e.g.][]{Khoshnevisan:2017, Khoshnevisan:2018, Kim:2019,Yi:2021}. Quite recently, our first main result (Theorem~\ref{THM:MAIN}) was used in~\cite{Yi:2021} to establish this property for certain versions of the stochastic heat and wave equations.  Specifically, the lower-tail bound in our Theorem~\ref{THM:MAIN} was sharp enough to enable exact calculations of the macroscopic Hausdorff dimension of exceedence sets of certain Gaussian random fields. (See Theorems 1.1 and 1.3, as well as Section 3.2 in~\cite{Yi:2021}.)

\subsection{Notation}
In order to simplify presentation, we always implicitly assume that $N\geq 2$, and we allow symbols for constants such as $c,C,c_0,C_0,\dots$ to change values with each appearance. When dealing with iterated logarithms, we use the abbreviation $\ell_2(N):=\max\{1,\log\log(N)\}$. For a real matrix $A$, define $\|A\|_F=\sqrt{\tr(A\ttop A)}$, $\|A\|_1=\sum_{i,j}|A_{ij}|$, and also define $\|A\|\op$ to be the largest singular value of $A$.
For a real vector $x$, we write $\|x\|_2$ for the Euclidean norm. If $A\in\R^{N\times N}$ is symmetric, its sorted eigenvalues are  $\lambda_{\max}(A)=\lambda_1(A)\geq\cdots\geq\lambda_N(A)=\lambda_{\min}(A)$. The identity matrix of size $N\times N$ is $I_N$, and the standard basis vectors in $\R^N$ are $e_1,\dots,e_N$. For the distribution of a random variable $V$, we write $\mathcal{L}(V)$, and we define its $\psi_1$-Orlicz norm as $\| V \|_{\psi_1} = \inf \{ t > 0 \, | \, \E [\exp (|V|/t)] \leq 2 \}$. If $\zeta$ is a Gaussian random vector with mean 0 and covariance matrix $\Sigma$, we write $\zeta\sim \mathcal{N}(0,\Sigma)$. For the univariate standard Gaussian distribution, the symbols $\Phi$ and $\phi$ denote the distribution function and density. Lastly, if $a$ and $b$ are real numbers, we use the notation $a\vee b=\max\{a,b\}$ and $a\wedge b=\min\{a,b\}$.

\section{Lower-tail bounds}\label{sec:main}

To clarify the statements of our main results, we first state a basic proposition describing the sizes of $\E(M_N(\xi))$ and $\text{med}(M_N(\xi))$ under condition~\eqref{eqn:assumcor}.  This proposition shows that the value $\sqrt{2(1-\rho_0)\log(N)}$ is a natural reference level for a lower-tail bound. Although this fact might be considered well-known among specialists, it is not easily referenced in the form given below, and so we provide a short proof at the end of Section~\ref{sec:mainproof}.

\begin{proposition}\label{prop:expec}
Let $\mu_N$ stand for either $\E(M_N(\xi))$ or $\textup{med}(M_N(\xi))$. If the condition~\eqref{eqn:assumcor} holds for some $\rho_0\in(0,1)$, then there is a universal constant $c_0>0$ such that 
\begin{equation}\label{eqn:expeclower}
\begin{split}
\ \ \ \ \ \ \ \ \ \ \ \ \ \ \ \mu_N & \ \geq \  \sqrt{2(1-\rho_0)\log(N)} - c_0\sqrt{\ell_2(N)}.
\end{split}
\end{equation}
Furthermore, if $R_{ij}=\rho_0$ for all $i\neq j$, then 
\begin{equation}\label{eqn:expecupper}
 \mu_N \ \leq \ \sqrt{2(1-\rho_0)\log(N)}.
\end{equation}
\end{proposition}

\noindent The following is our first main result, which will be extended and refined later in Corollary~\ref{cor:simple} and Theorem~\ref{thm:key}. The proof is deferred to Section~\ref{sec:mainproof}.
\begin{theorem}\label{THM:MAIN}
Fix two constants $\delta_0,\rho_0\in(0,1)$ with respect to $N$, and suppose the condition~\eqref{eqn:assumcor} holds. 
 Then, there is a constant $C>0$ depending only on $(\delta_0,\rho_0)$, such that
\begin{equation}\label{eqn:mainbound}
\small
  \P\bigg(M_N(\xi)\, \leq \,  \delta_0\sqrt{2(1-\rho_0)\log(N)}\bigg) \ \leq \  C\, N^{\ts\frac{-(1-\rho_0)(1-\delta_0)^2}{\rho_0}}(\log(N))^{\ts\frac{1-\rho_0(2-\delta_0)-\delta_0}{2\rho_0}}.
\end{equation}
Furthermore, the bound~\eqref{eqn:mainbound} is sharp in the sense that if $R_{ij}=\rho_0$ for all $i\neq j$, then there is a constant $c>0$ depending only on $(\delta_0,\rho_0)$, such that
 \begin{equation}\label{eqn:mainlower}
 \small
  \P\bigg(M_N(\xi) \, \leq \,  \delta_0\sqrt{2(1-\rho_0)\log(N)}\bigg) \ \geq \   c\, N^{\ts\frac{-(1-\rho_0)(1-\delta_0)^2}{\rho_0}}(\log(N))^{\ts\frac{1-\rho_0(2-\delta_0)-\delta_0}{2\rho_0}}.
\end{equation}
\end{theorem}
\noindent\textbf{Remarks.} To comment on some basic features of the theorem, first note that the dominant exponent   $-(1-\rho_0)(1-\delta_0)^2/\rho_0$ takes larger negative values as $\rho_0$ becomes smaller. 
Hence, the bound respects the fact that the lower-tail probability decays faster than any power of $N^{-1}$ when the entries of $\xi$ are independent. Also, the theorem conforms with the reference level motivated by Proposition~\ref{prop:expec}, since~\eqref{eqn:mainbound} implies that $\text{med}(M_N(\xi))$ cannot be much less than $\sqrt{2(1-\rho_0)\log(N)}$.

Regarding other recent lower-tail bounds for $M_N(\xi)$, their relation to~\eqref{eqn:mainbound} and~\eqref{eqn:mainlower} can be summarized as follows (in the setting of Theorem~\ref{THM:MAIN} with $\rho_0$ and $\delta_0$ regarded as fixed with respect to $N$). First, the paper~\citep[][Theorem 1.6]{Ding:2015} gives a two-sided tail bound for $M_N(\xi)$, implying that if $a$ and $b$ are positive constants satisfying  $\E(M_N(\xi))\geq a\sqrt{\log(N)}$ and $b\leq a/100$, then the probability $\P(M_N(\xi)\leq \E(M_N(\xi))-b\sqrt{\log(N)})$ is at most of order $N^{-b^2/(2-c(a))}$, for an unspecified constant $c(a)>0$.
Second, the paper~\citep[][Proposition 7]{Tanguy:2015} gives a two-sided tail bound implying that for any constant $b>0$, the probability $\P(M_N(\xi)\leq \E(M_N(\xi))-b\sqrt{\log(N)})$ is at most of order $\exp(-\ts\frac{cb}{\sqrt{\rho_0}}\sqrt{\log(N)})$, where $c>0$ is an unspecified constant. (Note that even if $c$ were specified, this bound would be of larger order than $N^{-\e}$ for any $\e\in(0,1)$.) Third, the paper~\citep[][Theorem 3.4]{Hartigan:2014} develops a lower-tail bound of the form $\P(M_N(\xi)\leq t_N(\alpha))\leq 2\alpha $ for any $\alpha\in (0,1/2)$, where $t_N(\alpha)$ is a threshold with a complex dependence on $N$, $\alpha$, as well as other parameters related to the correlation structure of $\xi$. In particular, Section 4 of the paper explains that the threshold can be expressed in terms of the minimum eigenvalue of $R$. For instance, when $R$ satisfies $R_{ij}=\rho_0$ for all $i\neq j$, the minimum eigenvalue is equal to $1-\rho_0$, and in this case, Theorem 3.4 and the comments preceding equation (4.1) in~\citep{Hartigan:2014} yield the following formula:
$t_N(\alpha)=\sqrt{1-\rho_0}[\kappa_N(\alpha)-\sqrt{\log(\kappa_N(\alpha))}+\sqrt{\rho_0}\Phi^{-1}(\alpha)]$, where we put $\kappa_N(\alpha)=\log(N^2/(2\pi))-2\log(-\log(\alpha))$.
 However, it seems that this formula for $t_N(\alpha)$ is not quite correct, since it is missing a square root. (Note that the quantity $\kappa_N(\alpha)$ scales like $\log(N)$, rather than $\sqrt{\log(N)}$, as a function of $N$.)\footnote{Relatedly, it also seems that the expression `$\lambda_n\sqrt{2\log(n)}$' in the paper's abstract should be replaced with $\sqrt{2\lambda_n\log(n)}$.}  Apart from this issue, the intricate form of the threshold also makes our result ostensibly easier to use. Lastly, in comparison to our work, none of the three mentioned papers resolve the question of whether or not the lower-tail bounds yield sharp rates.\\

\noindent\textbf{Extra  correlation structure.} We now present a direct corollary of Theorem~\ref{THM:MAIN} that allows for extra structure in the matrix $R$ to be used.
If the matrix $R$ satisfies $\max_{i\neq j}R_{ij}= \rho_0$, but ``most'' off-diagonal entries are substantially less than $\rho_0$, then we can gain considerable improvement.  
Roughly speaking, if there is a number $\rho_1\in (0,\rho_0)$ such that a large number of off-diagonal entries satisfy $R_{ij}\leq \rho_1$, then Theorem~\ref{THM:MAIN} can be improved by effectively replacing $\rho_0$ with the better value $\rho_1$.

\begin{corollary}\label{cor:simple}
Fix two constants $\delta_1,\rho_1\in(0,1)$ with respect to $N$. Suppose there is an index set $J\subset\{1,\dots,N\}$ with cardinality $|J|\geq 2$, such for any distinct $i,j\in J$, the bound $R_{ij}\leq \rho_1$ holds.
Then, there is a constant $ C>0$ depending only on $(\delta_1,\rho_1)$ such that
\small
\begin{equation*}
\P\bigg(M_N(\xi)\, \leq \,  \delta_1\sqrt{2(1-\rho_1)\log(|J|)}\bigg) \ \leq \ C\, |J|^{\ts\frac{-(1-\rho_1)(1-\delta_1)^2}{\rho_1} }\,(\log(|J|))^{\ts\frac{1-\rho_1(2-\delta_1)-\delta_1}{2\rho_1}}.
\end{equation*}
\normalsize
\end{corollary}
\noindent  Another point worth noting about this corollary is that it remains applicable even in situations where some variables are perfectly correlated and $\rho_0=1$.

\subsection{Further improvements using the restricted invertibility principle}\label{sec:ri}

In order to leverage the full strength of Corollary~\ref{cor:simple}, we need an index set $J\subset\{1,\dots,N\}$ with large cardinality, such that the off-diagonal entries $R_{ij}$ with $i,j\in J$ are uniformly small.  Quite remarkably, it turns out that such an index set is guaranteed to exist under rather general conditions,  as a consequence of a seminal result known as the \emph{restricted invertibility principle}~\citep[]{Bourgain:1987}.  
Over the past decade, this result has received considerable attention, due to its close links with the solution of the longstanding Kadison-Singer problem by Marcus, Spielman and Srivastava~\citep{Spielman:2015}. Below, we will apply a recent refinement of the restricted invertibility principle, established in a companion paper by the same group \citep{Spielman:2021}.

To introduce some notation, define the \emph{stable rank} ${\tt{r}}(A)$ of a non-zero positive semidefinite matrix $A$ as
$${{\tt{r}}}(A)= \frac{[\tr(A)]^2}{\|A\|_F^2}.$$
 This quantity always satisfies 
$$1\leq {\tt{r}}(A)\leq \textup{rank}(A),$$
and has the useful property that it approximately counts the number of dominant eigenvalues of $A$. In terms of the notion of stable rank, the restricted invertibility principle roughly says that for any matrix $L\in\R^{N\times N}$, there exists an index set $J\subset\{1,\dots,N\}$ with cardinality $|J|\approx {\tt{r}}(L\ttop L)$, such that 
the column submatrix of $L$ indexed by $J$ has singular values that are well separated from zero.

\begin{theorem}[Restricted invertibility principle~\text{\cite{Spielman:2021}}]\label{thm:ri} Let $L\in\R^{N\times N}$ be a non-zero matrix. Then, for any positive integer $k\leq {{\tt{r}}}(L\ttop L)$, there exists an index set $J\subset\{1,\dots,N\}$ with cardinality $|J|= k$,
 such that the following inequality holds for any real numbers $(a_j)_{j\in J}$,
\begin{equation}\label{eqn:ri}
\bigg\|\sum_{j\in J} a_j L e_j\bigg\|_2^2 \ \geq \ \frac{\|L\|_F^2}{N}\cdot\Big(1-\sqrt{\ts\frac{k}{{{\tt{r}}}(L\ttop L)}}\Big)^2\cdot\sum_{j\in J} a_j^2.
\end{equation}
\end{theorem}
We will apply the restricted invertibility principle above with $R^{1/2}$ playing the role of $L$. Note that because $R$ is a correlation matrix, we have $\|R^{1/2}\|_F^2/N=1$ and ${{\tt{r}}}(R)=N^2/\|R\|_F^2$. To simplify the inequality~\eqref{eqn:ri}, let $R(J)\in \R^{|J|\times |J|}$ denote the submatrix of $R$ with entries indexed by $J\times J$. Also, suppose there is an integer $k\geq 2$ and a scalar $\e\in (0,1)$, such that $k\leq \ts\frac{\e^2}{4}{{\tt{r}}}(R)$. In this case, the restricted invertibility principle ensures there is a choice of $J$ with cardinality equal to $k$ such that 
$$\lambda_{\min}(R(J)) \ \geq \ (1-\e/2)^2.$$
In turn, this implies that the off-diagonal entries of $R(J)$ are uniformly small, i.e.
\begin{equation*}
\begin{split}
\max_{i\neq j}R_{ij}(J) & \ = \  1-\min_{i\neq j}\ts\frac{1}{2}(e_i-e_j)\ttop R(J)(e_i-e_j)\\[0.2cm]
& \ \leq 1-\lambda_{\min}(R(J))\\[0.2cm]
& \ \leq \ \e.
\end{split}
\end{equation*}
The next theorem combines this information about $J$ with Corollary~\ref{cor:simple}. Later on, this connection will be used in our application dealing with rates of bootstrap approximation in high dimensions.

\begin{theorem}\label{thm:key}
Fix two constants $\delta,\e \in (0,1)$ with respect to $N$, and consider a Gaussian random vector \smash{$\xi\sim \mathcal{N}(0,R)$} for some correlation matrix  $R\in\R^{N\times N}$. Then, there is a constant $C>0$ depending only on $(\delta,\e)$ such that the following inequality holds for any integer $k\geq 2$ satisfying $k\leq \ts\frac{\e^2}{4} {{\tt{r}}}(R)$,
\begin{equation}\label{eqn:mainrefined}
\P\bigg(M_N(\xi)\, \leq \,  \delta \sqrt{2(1-\e)\log(k)}\bigg) \ \leq \ C\, k^{\frac{-(1-\e)(1-\delta)^2}{\e}}(\log(k))^{\frac{1-\e(2-\delta)-\delta}{2\e}}.
\end{equation}
\end{theorem}
\noindent\textbf{Remarks.} To discuss the choice of $k$, consider a basic situation where the root mean square of eigenvalues, say $\lambda_{\text{rms}}(R)=\ts\frac{1}{\sqrt N}\sqrt{\lambda_1(R)^2+\cdots+\lambda_N(R)^2}$, is upper bounded by a constant $C$. For instance, this condition is quite natural in the context of principal components analysis, where the matrix $R$ may only have a handful of dominant eigenvalues. Note too that this condition is much weaker than requiring $\|R\|\op\leq C$, since it is possible for the largest eigenvalue of $R$ to diverge while $\lambda_{\text{rms}}(R)$ remains bounded as $N\to\infty$. Due to the fact that $R$ is a correlation matrix, the condition $\lambda_{\text{rms}}(R)\leq C$ implies
$${\tt{r}}(R) \ \geq \ \ts\frac{N}{C^2},$$
and consequently, the integer $k$ in Theorem~\ref{thm:key} may be chosen to be of order $N$. This leads to an upper bound in~\eqref{eqn:mainrefined}
of order $N^{-(1-\e)(1-\delta)^2/\e}$, up to a logarithmic factor.
To see why this illustrates the benefit of the restricted invertibility principle, note that it enables us to bypass an assumption of the form $\max_{i\neq j}R_{ij}\leq \e$, which would have been necessary if Theorem \ref{THM:MAIN} had been used directly. Indeed, a condition like $\lambda_{\text{rms}}(R)\leq C$ is often more appealing from a statistical standpoint, because it allows for many off-diagonal entries of $R$ to be close to 1.\\

\section{Application to Bootstrap Methods in High Dimensions}\label{sec:app}

\noindent\textbf{Preliminaries.} Let $X_1,\dots,X_n\in\R^p$ be centered i.i.d.~observations with a standardized sum denoted as $S_n=n^{-1/2}\sum_{i=1}^n X_i$. In addition, define the coordinate-wise variances $\sigma_j^2=\var(X_{1j})$ for each $j=1,\dots,p$, which are assumed to be positive, 
 and define the partially standardized max statistic
\begin{equation*}
M=\max_{1\leq j\leq p} S_{n,j}/\sigma_j^{\tau},
\end{equation*}
where $\tau\in[0,1)$ is a tuning parameter. This type of max statistic was studied previously in~\citep{Lopes:2020}, and reduces to the basic example~\eqref{eqn:maxstatdef} in the case when $\tau$ is chosen to be 0.

In order to analyze bootstrap approximation of~$\mathcal{L}(M)$,  we will focus on the ``Gaussian multiplier bootstrap method''  popularized in~\citep{CCK:multiplier:2013}. To describe the method, define the sample mean $\bar X=\frac{1}{n}\sum_{i=1}^n X_i$, the sample covariance matrix $\hat\Sigma=\frac{1}{n}\sum_{i=1}^n (X_i-\bar X)(X_i-\bar X)\ttop$, and the coordinate-wise sample variances $\hat\sigma_j^2=\hat \Sigma_{jj}$ for $j=1,\dots,p$. The bootstrapped statistic $M^{\star}$ is constructed by generating a Gaussian random vector $S_n^{\star}\sim \mathcal{N}(0,\hat\Sigma)$, and defining\footnote{The expression $S_{n,j}^{\star}/\hat\sigma_j^{\tau}$ is regarded as 0 in the exceptional case $\hat\sigma_j=0$.}
\begin{equation*}
M^{\star}=\max_{1\leq j\leq p} S_{n,j}^{\star}/\hat\sigma_j^{\tau}.
\end{equation*}
Here, the key point is that it is possible in practice to directly simulate the conditional distribution of $M^{\star}$ given the data, denoted $\mathcal{L}(M^{\star}|X)$.  Accordingly, the multiplier bootstrap method uses $\mathcal{L}(M^{\star}|X)$ as an approximation to $\mathcal{L}(M)$.

As one last preliminary item, we will adopt the standard convention in high-dimensional statistics of considering a sequence of models indexed by $n$, in which all parameters may depend on $n$, except when stated otherwise. In particular, the dimension $p=p(n)$ is allowed to have arbitrary dependence on $n$. Likewise, if a parameter does not depend on $n$, then it does not depend on $p$ either. To simplify notation, we will write $a_n\lesssim b_n$ for two sequences of real numbers $a_n$ and $b_n$ if there is a constant $c>0$ not depending on $n$ such that $a_n\leq c b_n$ holds for all large $n$. In the case when both $a_n\lesssim b_n$ and $b_n\lesssim a_n$ hold, we write $a_n\asymp b_n$.

\subsection{A brief review of bootstrap approximation under variance decay}\label{sec:bootreview} Recently, the paper~\citep{Lopes:2020} analyzed how well $\mathcal{L}(M^{\star}|X)$ approximates $\mathcal{L}(M)$ in the setting of ``variance decay'', where the sorted coordinate-wise variances $\sigma_{(1)}^2\geq \cdots\geq\sigma_{(p)}^2$ have a polynomial decay profile. This means that there is a constant $\gamma>0$ not depending on $n$ such that the condition $\sigma_{(j)}\asymp j^{-\gamma}$ holds for all $j=1,\dots,p$. (The constant $\gamma$ is allowed to be arbitrarily large or small.) For instance, the structure of variance decay arises naturally in a variety of high-dimensional statistical applications related to \emph{principal components analysis}, \emph{count data}, and the \emph{Fourier coefficients of functional data}. Under the assumption of variance decay, as well as some additional assumptions on the correlation and tail-behavior of the covariates, Theorem 3.2 in the paper~\citep{Lopes:2020} established the following bootstrap approximation result. Namely, for any fixed $\delta\in (0,1/2)$, there is a constant $C>0$ not depending on $n$ such that the bound

\begin{equation}\label{eqn:llm}
\sup_{t\in\R}\Big|\P(M^{\star}\leq t|X)-\P(M\leq t)\Big| \ \leq \ C n^{-\frac{1}{2}+\delta}
\end{equation}
holds with probability at least $1-C/n$. 

The result~\eqref{eqn:llm} has some significant distinguishing features in relation to the other recent works~\citep{Lopes2:2020,Chernozhukov:2020} that have obtained near $n^{-1/2}$ rates of bootstrap approximation for max statistics in high-dimensional settings. First, the use of variance decay makes it possible for the bound~\eqref{eqn:llm} to hold independently of $p$, whereas the other works do not use this structure and obtain rates with polylogarithmic dependence on $p$.
 Second, the other works require that the correlation matrix $\text{cor}(X_1)\in\R^{p\times p}$ be either positive definite or well-approximated by a positive definite matrix, whereas the bound~\eqref{eqn:llm} does not. Nevertheless, the bound~\eqref{eqn:llm} does rely on some other assumptions about the matrix $\text{cor}(X_1)$. Below, we outline these assumptions, and later in Section~\ref{sec:newboot}, we will show how these assumptions can be simplified and weakened by applying our work from Section~\ref{sec:main}.
To introduce some notation, consider an arbitrary index $d\in\{1,\dots,p\}$, and let $J(d)\subset \{1,\dots,p\}$ denote a set of indices corresponding to the $d$ largest coordinate-wise variances, i.e. $\{\sigma_{(1)}^2,\dots,\sigma_{(d)}^2\}=\{\sigma_j^2\,|\,j\in J(d)\}$. Next, let $R(d)$ denote the $d\times d$ correlation matrix of the variables $(X_{1j})_{j\in J(d)}$, and let the matrix $R^+(d)\in\R^{d\times d}$ be the non-negative part, $R_{ij}^+(d)=R_{ij}(d)\vee 0$. Lastly, define the integer $m_n=\lceil (\log(n))^3\wedge p\rceil$, which always satisfies $1\leq m_n \leq p$. 

In terms of this notation, the paper~\citep{Lopes:2020} makes the following three correlation assumptions in order to establish~\eqref{eqn:llm}:
\begin{enumerate}[(a)]
\item There is a constant $\e_0\in (0,1)$ not depending on $n$ such that 
$$\max_{i\neq j} R_{ij}(m_n)\leq 1-\e_0.$$
\vspace{-0.3cm}
\item The matrix $R^+(m_n)$ is positive semidefinite.\\[-0.cm]
\item The condition $\|R(m_n)\|_1\lesssim m_n$ holds.
\end{enumerate}
Concerning assumption (b), this is non-trivial because the operation of taking the entrywise non-negative part of a matrix does not generally preserve positive semidefiniteness~\citep[][Theorem 4.11]{Guillot:2015}.

\subsection{Bootstrap approximation result with relaxed correlation assumptions}\label{sec:newboot}
 In this subsection, we will present a result of the form~\eqref{eqn:llm} in which the previous correlation assumptions (a) and (b) are removed, and in which (c) will be replaced by a version involving a weaker norm. Apart from the correlation structure, the following conditions on the data-generating model and the variance decay profile are slightly simplified versions of the ones used in~\citep{Lopes:2020}.

 \begin{assumption}[Data-generating model and variance decay]\label{assump1}~ \vspace{0.1cm} %
 \begin{enumerate}[(i)]
\item  There is a positive semidefinite matrix $\Sigma \in \R^{p \times p}$, such that the observations $X_1, \dots, X_n \in \R^p$ are generated as 
$$X_i = \Sigma^{1/2} Z_i$$
 for each $i=1,\dots,n$, where $Z_1, \dots, Z_n$ are i.i.d.~random vectors, with $\E (Z_1) = 0$, $\E(Z_1 Z_1^\top) = I_p$, as well as $\sup_{ \| u\|_2 = 1} \| Z_1^\top u\|_{\psi_1} \lesssim 1$.\\[-0.2cm]
\item  The parameters $\sigma_{(1)}\geq \cdots\geq \sigma_{(p)}$ are positive, and there is a constant $\gamma>0$ not depending on $n$ such that the condition
\[  \sigma_{(j)} \asymp  j^{- \gamma} \]
holds for all $j=1,\dots,p$. 
\end{enumerate}
\end{assumption}

The following is the main result of the current section, and the proof will be given in Section~\ref{sec:bootproof}.

 \begin{theorem}[Bootstrap approximation] \label{THM:BOOT}
Fix any constants $\delta\in (0,1/2)$ and $\tau\in [0,1)$ with respect to $n$, and define  $\kappa=4\vee (3\gamma(1-\tau))$ as well as $l_n= \big\lceil n^{\delta/\kappa}\wedge p\big\rceil$. In addition, suppose that Assumption~\ref{assump1} holds, and that the condition
\begin{equation}\label{eqn:newcor}
\|R(l_n)\|_F^2 \ \lesssim \ l_n^{2-\delta}.
\end{equation}
holds. Then, there is a constant $C>0$ not depending on $n$ such that the event
\begin{equation}\label{eqn:newllm}
\sup_{t\in\R}\Big|\P(M^{\star}\leq t|X)-\P(M\leq t)\Big| \ \leq \ Cn^{-\frac{1}{2}+\delta}
\end{equation}
holds with probability at least $1-C/n$.
\end{theorem}

\noindent\textbf{Remarks.} To interpret the correlation assumption~\eqref{eqn:newcor}, it should be emphasized that the upper bound $\|R(l_n)\|_F^2\leq l_n^2$ holds under all circumstances, because any positive semidefinite matrix $A$ satisfies $\|A\|_F^2\leq [\tr(A)]^2$. So, given that $\delta$ may be taken to be arbitrarily small, it is not possible to substantially weaken~\eqref{eqn:newcor} in general.  Furthermore, it is worth noting that the condition~\eqref{eqn:newcor} only applies to the small set of variables indexed by $J(l_n)$, while all other variables indexed by $\{1,\dots,p\}\setminus J(l_n)$ have a correlation structure that is unrestricted. 

A large class of examples of correlation matrices satisfying~\eqref{eqn:newcor} can be constructed in the following way. Let $f:[0,\infty)\to [0,1]$ be any continuous convex function that satisfies $f(0)=1$, as well as $f(t)\leq  ct^{-\delta}$ for some fixed constants $c,\delta>0$ and all $t\geq 0$. Then, by P\'olya's criterion~\citep[][Theorem 1]{Polya}, any matrix with entries defined by $R_{ij}=f(|i-j|)$ will be a correlation  matrix that satisfies~\eqref{eqn:newcor}. Moreover, any other correlation matrix obtained by permuting the rows and columns will continue to satisfy~\eqref{eqn:newcor}.  For additional discussion of such matrices, and some of their connections to high-dimensional statistics, we refer to~\citep{Bickel:2008:a}.

With regard to our lower-tail bounds from Section~\ref{sec:main}, the connection with Theorem~\ref{THM:BOOT} can be explained as follows. Overall, the proof of this result is based on a dimension-reduction strategy that involves showing $M$ is well approximated by a statistic of the form $M'=\max_{j\in J'}S_{n,j}/\sigma_j^{\tau}$, where $J'\subset\{1,\dots,p\}$ is an index set with cardinality much less than $p$. In order to show that $M=M'$ holds with high probability,  we will localize the maximizing index for $M$ to the set $J'$. That is, if $\hat{\j}$ is a random index such that $M=S_{n,\hat{\j}}/\sigma_{\hat{\j}}^{\tau}$, then we will show that $\hat{\j}$ falls into $J'$ with high probability. It is in this localization step that the lower-tail bound from Theorem~\ref{thm:key} will be used: Two values $x<y$ will be determined, such that $M'$ is likely to be above $y$, while the maximum of $S_{n,j}/\sigma_j^{\tau}$ over the indices $j\not\in J'$ is likely to be below $x$. Hence, the first of these items requires a lower-tail bound, because it involves showing that the probability $\P(M'\leq y)$ is small.

\section{Proof of Theorem~\text{\ref{THM:MAIN}}}\label{sec:mainproof}

Throughout the proofs, we always assume that $N$ is sufficiently large for any given expression to make sense --- since this can be accommodated by an adjustment of the constants $C$ and $c$ in the statement of Theorem~\ref{THM:MAIN}.  The symbols $c$ and $C$ will also be frequently be re-used with the understanding that they never depend on $N$ (and similarly with respect to the sample size $n$ in Section~\ref{sec:bootproof}). In addition, it will simplify presentation to introduce the the quantities
\begin{equation*}
\alpha_0=\ts\frac{(1-\rho_0)(1-\delta_0)^2}{\rho_0} \ \ \ \ \ \text{ and }\  \ \ \ \  \beta_0 = \ts\frac{(1-\rho_0)(1-\delta_0)}{\rho_0},
\end{equation*}
so that the bounds in Theorem~\ref{THM:MAIN} are proportional to $N^{-\alpha_0}(\log(N))^{\frac{\beta_0-1}{2}}$.\\[-0.2cm]

\noindent\textbf{Proof of the upper bound~\eqref{eqn:mainbound}.} Let $t_N= \delta_0\sqrt{2(1-\rho_0)\log(N)}$, and let $\xi'\sim \mathcal{N}(0,R')$  be a Gaussian vector in $\R^N$, where $R'$ is a correlation matrix satisfying $R_{ij}'=\rho_0$ for all $i\neq j$. Due to the assumption that $\max_{i\neq j}R_{ij}\leq \rho_0$, it follows from Slepian's Lemma~\citep{Slepian:1962} that 
\begin{equation*}
 \P(M_N(\xi)\leq t_N)  \ \leq \ \P(M_N(\xi')\leq t_N).
\end{equation*}
To control the larger probability, let $\zeta_0,\zeta_1,\dots,\zeta_N$ denote independent $\mathcal{N}(0,1)$ random variables, and note that the coordinates of $\xi'$ may be jointly represented in distribution as
\begin{equation}\label{eqn:correp}
\xi_j' = \sqrt{\rho_0}\zeta_0+\sqrt{1-\rho_0}\zeta_j.
\end{equation}
This yields the following representation of the maximum
\begin{equation}\label{eqn:maxrep}
M_N(\xi')= \sqrt{\rho_0}\zeta_0+\sqrt{1-\rho_0}\max_{1\leq j\leq N}\zeta_j,
\end{equation}
which allows us to view $M_N(\xi')$ as the convolution of the two independent random variables on the right side. Consequently, a direct calculation may be used to obtain the exact formula
\begin{equation}\label{eqn:mainint}
\P\big(M_N(\xi') \leq t_N\big) = \int_{-\infty}^{\infty}\psi_N(s)ds,
\end{equation}
where the integrand is defined by
$$\psi_N(s) \ = \ \phi(s)\Phi^N \big(\ts\frac{t_N-\sqrt{\rho_0}s}{\sqrt{1-\rho_0}}\big),$$
with $\phi$ and $\Phi$ being the standard normal density and distribution function.
The remainder of the proof consists in bounding integral $\int_{-\infty}^{\infty}\psi_N(s)ds$ over the intervals 
$(-\infty,-c_N]$, $[-c_N,0]$, and $[0,\infty)$, where we define
\begin{equation}\label{eqn:cnformula}
c_N=\sqrt{\ts\frac{2(1-\rho_0)\log(N)}{\rho_0}}\Big\{1-\delta_0-\ts\frac{\frac 14 \ell_2(N)}{\log(N)}\Big\}.
\end{equation}

\noindent\emph{Remarks.} Our choice of the cut-off point $c_N$ is a crucial element of the proof. To give a sense of how delicate this is, a close inspection of the proof shows that the $1/4$ coefficient on the trailing term $\ell_2(N)/\log(N)$ is needed to obtain both the upper and lower bounds in Theorem~2.1. Some intuition for the definition of $c_N$ can be given as follows. First, define $b_N=\sqrt{2\log(N)}-\ell_2(N)/\sqrt{8\log(N)}$ and the random variable
\smash{
$G_N=\sqrt{2\log(N)}(\max_{1\leq j\leq N}\zeta_j-b_N),$}
which has the property that $G_N+\log(\sqrt{4\pi})$ converges weakly to a standard Gumbel distribution as $N\to\infty$~\citep[][Theorem 1.5.3]{Leadbetter:1983}.
Based on this notation and~(4.2), it can be checked that 
\begin{equation}\label{eqn:heuristic}
\P(M_N(\xi')\leq t_N) \ =  \  \P\Big(\zeta_0 + \sqrt{\ts\frac{1-\rho_0}{2\rho_0\log(N)}}G_N \leq -c_N\Big).
\end{equation}
Next, the probability on the right can be approximated heuristically by dropping the random variable $\sqrt{(1-\rho_0)/(2\rho_0)}\,G_N/\sqrt{\log(N)}$, because it is of order $1/\sqrt{\log(N)}$ (in probability), and hence asymptotically negligible compared to $\zeta_0$. After this heuristic is used, some further calculation then leads to the surmise that $\P(M_N(\xi')\leq t_N)$ should be of order $N^{-\alpha_0}\log(N)^{\frac{\beta_0-1}{2}}$. However, this line of reasoning will not be used in the formal proof, because the bounds in Theorem~2.1 are non-asymptotic. Also, note that this asymptotic heuristic suppresses the fact that $1/\sqrt{\log(N)}$ approaches 0 very slowly, which matters much more from a non-asymptotic standpoint.

Returning to the main thread, the problem of upper bounding the integral of $\psi_N(s)$ over the interval $[-c_N,0]$ is the most involved part of the proof, and is postponed to Lemma~\ref{lem:hard} later on. (The reason for this difficulty is that the function $\psi_N(s)$ is especially sensitive to small changes in $s$ over $[-c_N,0]$. In particular, this stage of the analysis involves developing separate bounds over two sub-intervals of $[-c_N,0]$ in order to account for changes in the local behavior of the function.)
Once the proof of that lemma is complete, it will be straightforward to control the integral over $[0,\infty)$, which is done subsequently in Lemma~\ref{lem:easy}. For the moment, we only handle the interval $(-\infty,-c_N]$, since it requires no further preparation. Indeed, we have
\begin{equation}\label{eqn:firstmills}
\begin{split}
\int_{-\infty}^{-c_N} \psi_N(s) ds & \ \leq \ \int_{-\infty}^{-c_N}\phi(s)ds\\[0.1cm]
& \ \leq \ \frac{e^{-c_N^2/2}}{\sqrt{2\pi}c_N} \ \ \ \ \text{(Mill's inequality)}\\[0.3cm]
& \ \leq \ \ts\frac{C}{\sqrt{\log(N)}}\,  N^{-\alpha_0} \, \exp\big\{\frac{\beta_0}{2} \ell_2(N)\big\}\\[0.3cm]
& \ \leq \ C N^{-\alpha_0} (\log(N))^{\frac{\beta_0-1}{2}},
\end{split}
\end{equation}
where $C$ only depends on $\delta_0$ and $\rho_0$. This proves the upper bound~\eqref{eqn:mainbound}.\\

\noindent\textbf{Proof of the lower bound~\eqref{eqn:mainlower}.} It suffices to derive a lower bound on the integral $\int_{-\infty}^{-c_N}\psi_N(s)ds$. Using the fact that $\Phi^N \Big(\ts\frac{t_N-\sqrt{\rho_0}s}{\sqrt{1-\rho_0}}\Big)$ is decreasing in $s$, we have
\begin{equation}\label{eqn:firstlower}
\int_{-\infty}^{-c_N}\psi_N(s)ds \ \geq \ \Phi^N \Big(\ts\frac{t_N+\sqrt{\rho_0}c_N}{\sqrt{1-\rho_0}}\Big) \Phi(-c_N).
\end{equation}
If we define the number
\begin{equation}
\begin{split}
r_N& \ = \ \ts\frac{t_N+\sqrt{\rho_0}c_N}{\sqrt{1-\rho_0}} \\
& \ = \ \sqrt{2\log(N)}\Big\{1-\ts\frac{\frac 14\ell_2(N)}{\log(N)}\Big\},
\end{split}
\end{equation}
then Mill's inequality gives
\begin{equation}
\begin{split}
\Phi(r_N) &  \ \geq \ 1-\ts\frac{1}{\sqrt{2\pi}r_N} e^{-r_N^2/2}\\[0.2cm]
 & \ \geq   \ 1 -\ts\frac{C}{N}.
\end{split}
\end{equation}
Hence, the limit $(1-C/N)^N\to \exp(-C)$ as $N\to\infty$ shows that $\Phi^N(r_N)$ is lower bounded by a positive constant for all large $N$. Finally, regarding the factor $\Phi(-c_N)$ in~\eqref{eqn:firstlower}, the lower-bound version Mill's inequality~\cite[Theorem 1.4]{Durrett:2005} gives
\begin{equation}
\begin{split}
\Phi(-c_N) &  \ \geq \  \Big\{\ts\frac{1}{c_N}-\ts\frac{1}{c_N^3}\Big\}\frac{1}{\sqrt{2\pi}}e^{-\frac{1}{2}c_N^2}\\[0.2cm]
& \  \geq \ \ c N^{-\alpha_0}(\log(N))^{\frac 12(\beta_0-1)},
\end{split}
\end{equation}
where the second step follows from a calculation very similar to~\eqref{eqn:firstmills}. This completes the proof of the lower bound~\eqref{eqn:mainlower} and the theorem. \qed

\begin{lemma}\label{lem:hard}
There is a constant $C>0$ depending only on $(\rho_0,\delta_0)$ such that
\begin{equation}
\int_{-c_N}^{0}\psi_N(s)ds \ \leq \ C\, N^{-\alpha_0} (\log(N))^{\frac 12(\beta_0-1)} .
\end{equation}
\end{lemma}

\proof 
To simplify the analysis of $\psi_N(s)$, we will introduce a change of variable and write $s$ as function of a number $\delta$ lying in an interval denoted as
$$\mathcal{I}_N=\Big[0\, , \, 1-\delta_0-\ts\frac{\frac 14 \ell_2(N)}{\log(N)}\Big].$$
Specifically, we write
$s=s_N(\delta)  = -a_N\delta,$
where
$a_N=\sqrt{\ts\frac{2(1-\rho_0)\log(N)}{\rho_0}}.$
Also note that $s_N(\cdot)$ maps $\mathcal{I}_N$ to the interval $[-c_N,0]$. As a means of simplifying the factor $\Phi^N\big(\ts\ts\frac{t_N-\sqrt{\rho_0}s}{\sqrt{1-\rho_0}}\big)$ in the definition of $\psi_N(s)$, let
\begin{equation}\label{eqn:undef}
\begin{split}
u_N(\delta) & \ = \ \ts\frac{t_N-\sqrt{\rho_0}s_N(\delta)}{\sqrt{1-\rho_0}} \
= \ \sqrt{2\log(N)}(\delta_0+\delta).
\end{split}
\end{equation}
It follows that for all $\delta\in\mathcal{I}_N$,
\begin{equation}\label{eqn:logpsiformula}
\log \psi_N(s_N(\delta)) \ = \  -\ts\frac{(1-\rho_0)\log(N)\delta^2}{\rho_0}+N\log \Phi(u_N(\delta))-\log(\sqrt{2\pi}). 
\end{equation}
Due to the lower-bound form of Mill's inequality, we have
\begin{equation}\label{eqn:cubicmills}
\begin{split}
\Phi(u_N(\delta)) &\ \leq \ 1-\ts\frac{1}{\sqrt{2\pi}}\Big\{\ts\frac{1}{u_N(\delta)}-\ts\frac{1}{u_N^3(\delta)}\Big\}\displaystyle\exp({-\ts\frac{1}{2}u_N(\delta)^2})\\[0.2cm]
 & \ = \ 1-w_N(\delta),
\end{split}
\end{equation}
where $w_N(\delta)$ is defined by the last line. When $N$ is sufficiently large, it is simple to check that the condition $0<w_N(\delta)<1$ holds for all $\delta\in\mathcal{I}_N$, which gives
$\log(1-w_N(\delta))  \leq  -w_N(\delta).$
Combining the last few steps, the following bound holds for all $\delta\in\mathcal{I}_N$,
\begin{equation}\label{eqn:logupperbound}
\begin{split}
N\log \Phi(u_N(\delta)) \ & \ \leq -N\, w_N(\delta)\\[0.2cm]
& \ \leq \ \ts\frac{-cN}{\sqrt{\log(N)}}\exp\big\{-\log(N)(\delta_0+\delta)^2\big\}.
\end{split}
\end{equation}

The work up to this point provides us with a useful majorant for $\psi_N$. By looking at the equation~\eqref{eqn:logpsiformula} and the bound~\eqref{eqn:logupperbound}, it is clear that if we define the function
\begin{equation}\label{eqn:upperfunc}
\begin{split}
f_N(\delta)&=\ts\frac{(1-\rho_0)\log(N)\delta^2}{\rho_0}+\ts\frac{cN}{\sqrt{\log(N)}}\displaystyle\exp\big\{{-\log(N)(\delta_0+\delta)^2}\big\},
\end{split}
\end{equation}
then the bound
\begin{equation}\label{eqn:upperfuncbound}
\psi_N(s_N(\delta)) \ \leq \ \exp\{-f_N(\delta)\}
\end{equation}
holds for all $\delta\in\mathcal{I}_N$. Integrating this bound gives
\begin{equation}\label{eqn:firstcov}
\int_{-c_N}^0 \psi_N(s)ds \ \leq \ a_N \!\int_{\mathcal{I}_N} \exp\{-f_N(\delta)\} d\delta.
\end{equation}
We now introduce another change of variable, and write $\delta$ as a function of a positive number $\eta$ using  
$$\delta=\delta_N(\eta) = 1-\delta_0-\ts\frac{\eta\ell_2(N)}{\log(N)}.$$
If we define the interval
$\J_N =\big[\ts\frac{1}{4}\, , \, \frac{(1-\delta_0)\log(N)}{\ell_2(N)}\big],$
then $\delta_N(\cdot)$ maps $\J_N$ to $\mathcal{I}_N$, and the integral bound~\eqref{eqn:firstcov} becomes
\begin{equation}\label{eqn:secondcov}
\begin{split}
\int_{-c_N}^0 \psi_N(s)ds  
& \ \leq \ts\frac{a_N \ell_2(N)}{\log(N)} \displaystyle \int_{\J_N} \exp\{-f_N(\delta_N(\eta))\}d\eta.
\end{split}
\end{equation}

The remainder of the proof will be divided into two parts, in which the integral over $\J_N$ is decomposed with the subintervals 
$$\J_N' = \Big[\ts\frac{1}{4}\, , \, \ts\frac{(\log(N))^{1/4}}{\ell_2(N)}\Big] \ \ \ \ \ \text{ and } \ \ \ \ \ \J_N'' = \Big[\ts\frac{(\log(N))^{1/4}}{\ell_2(N)},\ts\frac{(1-\delta_0)\log(N)}{\ell_2(N)}\Big].$$
In handling these subintervals below, it will be convenient to label the summands of $f_N$ in line~\eqref{eqn:upperfunc} according to
$$f_N(\delta) = g_N(\delta) + h_N(\delta),$$
where
\begin{align*}
g_N(\delta) \ & \ = \ts\frac{(1-\rho_0)\log(N)\delta^2}{\rho_0}\\[0.2cm]
h_N(\delta) \ &\ =\ts\frac{cN}{\sqrt{\log(N)}}\displaystyle\exp\big\{{-\log(N)(\delta_0+\delta)^2}\big\}.
\end{align*}

\noindent\textbf{The integral over $\J_N'$.}
By expanding out the square $\delta_N(\eta)^2$, and dropping the smallest positive term, the following bound holds for any $\eta$,
$$\exp\big\{\!-\!g_N(\delta_N(\eta))\big\} \ \leq   \ N^{-\alpha_0}\exp\Big\{2\beta_0 \ell_2(N)\eta\Big\}.$$
In addition, if we expand the square $(\delta_0+\delta_N(\eta))^2$, and use the fact that every $\eta\in\J_N'$ is bounded above by $(\log(N))^{1/4}/\ell_2(N)$, then a short calculation gives
\begin{equation*}
\begin{split}
-h_N(\delta_N(\eta))
& \ \leq \ -c e^{\ell_2(N)(2\eta-\frac{1}{2})},
\end{split}
\end{equation*}
for small enough $c>0$. Directly combining the last two steps gives
\begin{equation}\label{eqn:quadstep}
\exp\big\{\!-\!f_N(\delta_N(\eta))\big\} \ \leq  \ N^{-\alpha_0}\exp\Big\{ 2\beta_0\eta \ell_2(N)  - ce^{\ell_2(N)(2\eta-\frac{1}{2})} \Big\}.
 \end{equation}
To simplify the previous bound, define
$x(\eta)=2\eta-\ts\frac{1}{2}.$
Since $x(\eta)$ is non-negative for all $\eta\in\J_N'$, we may approximate $\exp\{\ell_2(N)x(\eta)\}$ from below using a second-order Taylor expansion
$1+\ell_2(N) x(\eta)+\ts\frac{1}{2}\ell_2(N)^2 x(\eta)^2$.
After some arithmetic, the bound~\eqref{eqn:quadstep} becomes
\begin{equation}\label{eqn:expbound}
\begin{split}
\exp\big\{\!-\!f_N(\delta_N(\eta))\big\} & \ \leq \ e^{\frac{\beta_0}{2}\ell_2(N)}\cdot  e^{-c} \cdot N^{-\alpha_0}\cdot \varphi_N(x(\eta)),
 \end{split}
 \end{equation}
where we define the function
$$\varphi_N(x)=\exp\Big\{ -\big[\ts\frac{c}{2}\ell_2(N)^2\big] x^2 + \big[(\beta_0-c)\ell_2(N)\big] x  \Big\}.$$
Integrating the bound~\eqref{eqn:expbound} over $\J_N'$, we obtain
\begin{equation}\label{eqn:nearly}
\footnotesize
\begin{split}
 \frac{a_N \ell_2(N)}{\log(N)} \displaystyle \int_{\J_N'}\exp\{-f_N(\delta_N(\eta))\}d\eta  
&  \ \ \leq \ \ C\cdot \ell_2(N) \cdot (\log(N))^{\frac{1}{2}(\beta_0-1)} \cdot N^{-\alpha_0} \displaystyle\int_{0}^{\infty} \varphi_N(x)dx.
 \end{split}
 \end{equation}
To handle the last integral, note that the function $\varphi_N(x)$ can be written in the form $\varphi_N(x)=\exp(-a x^2+bx)$, 
and that the elementary Gaussian integral bound  
\begin{equation*}
\int_0^{\infty} e^{-a x^2+b x}dx \ \leq \  \ts\frac{\sqrt{\pi}}{\sqrt a} \exp\big(\ts\frac{b^2}{4a}\big)
\end{equation*}
holds for any $a>0$ and $b\in\R$.
Therefore,
$$ \int_{0}^{\infty}\varphi_N(x)dx \ \leq \ \ts\frac{C}{\ell_2(N)}. $$
Combining this with the bound~\eqref{eqn:nearly} completes the work on $\J_N'$.\\

\noindent\textbf{The integral over $\J_N''$.} First note that $\exp\{-f_N(\delta)\} \leq \exp\{-h_N(\delta)\}$. Also, the function $\exp\{-h_N(\delta_N(\eta))\}$ is decreasing in $\eta$, and so if we denote the left endpoint of $\J_N''$ as 
$$\eta_N=(\log(N))^{1/4}/\ell_2(N),$$ then we have the following height-width integral bound
\begin{equation*}
\begin{split}
\ts\frac{a_N \ell_2(N)}{\log(N)} \displaystyle \int_{\J_N''}\exp\big\{\!-\!f_N(\delta_N(\eta))\big\}d\eta & \ \leq \ \ts\frac{a_N \ell_2(N)}{\log(N)} \cdot |\J_N''| \cdot \exp\{-h_N(\delta_N(\eta_N))\}\\[0.3cm]
& \ \leq \ C\cdot \sqrt{\log(N)}  \cdot \exp\Big\{\ts\frac{-c}{\sqrt{\log(N)}}e^{2(\log(N))^{1/4}} \Big\}.
\end{split}
\end{equation*}
This bound is of smaller order than $N^{-\alpha_0}(\log(N))^{\frac 12(\beta_0-1)}$,  which completes the proof of Lemma~\ref{lem:hard}.\qed

\begin{lemma}\label{lem:easy}
There is a constant $C>0$ depending only on $(\delta_0,\rho_0)$ such that 
\begin{equation}
\int_{0}^{\infty}\psi_N(s)ds \ \leq \ C\, N^{-\alpha_0} (\log(N))^{\frac 12(\beta_0-1)} .
\end{equation}
\end{lemma}

\proof The function $\Phi^N \Big(\ts\frac{t_N-\sqrt{\rho_0}s}{\sqrt{1-\rho_0}}\Big)$ is decreasing in $s$, and so
$$\int_0^{\infty}\psi_N(s)ds \leq \  \Phi^N \Big(\ts\frac{t_N}{\sqrt{1-\rho_0}}\Big).$$
Then, the inequality~\eqref{eqn:logupperbound} gives
\begin{equation}
\begin{split}
 \Phi^N \Big(\ts\frac{t_N}{\sqrt{1-\rho_0}}\Big) & \ \leq \ \exp(-h_N(0))\\[0.2cm]
&  \ = \ \exp\Big\{\ts\frac{-c N^{1-\delta_0^2}}{\sqrt{\log(N)}}\Big\},
 \end{split}
 \end{equation}
 which is clearly of smaller order than the stated bound. \qed

~\\

\noindent\textbf{Proof of Proposition~\ref{prop:expec}.} We first show~\eqref{eqn:expeclower} in the case where $\mu_N=\E(M_N(\xi))$. Let $R'\in\R^{N\times N}$ be a correlation matrix with $R_{ij}'=\rho_0$ for all $i\neq j$, and let $\xi'\sim \mathcal{N}(0,R')$. Based on the assumed condition~\eqref{eqn:assumcor}, it follows from Slepian's lemma that
\begin{equation}
 \E(M_N(\xi))  \ \geq \  \E(M_N(\xi')).
\end{equation}
Next, the representation~\eqref{eqn:maxrep} gives
\begin{equation}\label{eqn:twostep}
\small
\begin{split}
\E\big[M_N(\xi')\big]  & \ = \ \sqrt{1-\rho_0}\, \E\Big[\max_{1\leq j\leq N}\zeta_j\Big]\\[0.3cm]
& \ \geq \  \sqrt{2(1-\rho_0)\log(N)} - c_0\sqrt{\ell_2(N)}
\end{split}
\end{equation}
for some universal constant $c_0>0$, where the second step can be obtained from~\cite[p.66]{Massart:2007}.
To handle $\mu_N=\text{med}(M_N(\xi))$, note that
\begin{equation}
\begin{split}
\Big|\E(M_N(\xi))-\text{med}(M_N(\xi))\Big| & \ \leq \ \sqrt{\var(M_N(\xi))}\\[0.2cm]
& \ \leq 1,
\end{split}
\end{equation}
where the second step uses the Gaussian Poincar\'e inequality~\cite[Theorem 3.20]{Massart:2013}.
To show the upper bound~\eqref{eqn:expecupper}, note that we may take $\xi=\xi'$. When $\mu_N=\E(M_N(\xi'))$, the result follows from the first part of~\eqref{eqn:twostep} and the standard inequality $\E(\max_{1\leq j\leq N}\zeta_j) \leq \sqrt{2\log(N)}$. When $\mu_N=\text{med}(M_N(\xi'))$, the result follows from the fact that $\text{med}(M_N(\xi'))\leq \E(M_N(\xi'))$, since $M_N(\cdot)$ is a continuous convex function~\citep{Kwapien:1994}. \qed

\section{Proof of Theorem~\ref{THM:BOOT}} \label{sec:bootproof}
Let $\Sigma=\E(X_1X_1\ttop)$ be the population covariance matrix, and let $\hat\Sigma$ be the sample covariance matrix defined in Section~\ref{sec:app}, with $\hat\sigma_j^2=\hat\Sigma_{jj}$ for every $j=1,\dots,p$. Recall that $S_n^{\star}$ denotes a Gaussian random vector drawn from $\mathcal{N}(0,\hat\Sigma)$, and define  a corresponding Gaussian vector $\tilde S_n\sim\mathcal{N}(0,\Sigma)$.  For each index $d\in\{1,\dots,p\}$, let the index set $J(d)$ be as defined in Section~\ref{sec:app} and define three associated max statistics%
\begin{align*}
M_{d}&=\max_{j\in J(d)} S_{n,j}/\sigma_j^{\tau}\\[0.2cm]
\tilde M_{d}&=\max_{j\in J(d)} \tilde S_{n,j}/\sigma_j^{\tau}\\[0.2cm]
M_{d}^{\star}&=\max_{j\in J(d)} S_{n,j}^{\star}/\hat\sigma_j^{\tau}.
\end{align*}
To compare these statistics, we will use the Kolmogorov metric, defined for generic random variables $U$ and $V$ as
$$d_{\textup{K}}(\mathcal{L}(U),\mathcal{L}(V)) = \sup_{t\in\R}\big|\P(U\leq t)-\P(V\leq t)\big|.$$

In this notation, the proof amounts to showing there is a constant $C>0$ not depending on $n$ such that the event 
$$d_{\textup{K}}(\mathcal{L}(M_{p}),\mathcal{L}(M_{p}^{\star}|X)) \ \leq \ Cn^{-\frac{1}{2}+\delta}$$
holds with probability at least $1-C/n$.
The overall structure of the proof will decomposed into two parts by considering the triangle inequality
$$d_{\textup{K}}(\mathcal{L}(M_{p}),\mathcal{L}(M_{p}^{\star}|X)) \ \leq \ d_{\textup{K}}(\mathcal{L}(M_{p}),\mathcal{L}(\tilde M_{p})) \ + \ d_{\textup{K}}(\mathcal{L}(\tilde M_{p}),\mathcal{L}(M_{p}^{\star}|X))
$$
and then separately bounding the two distances on the right. Hence, it suffices to prove the following proposition.

\begin{proposition}\label{PROP:BOOT}
Suppose the conditions of Theorem~\ref{THM:BOOT} hold. Then, 
\begin{equation}\label{eqn:gpart}
 d_{\textup{K}}(\mathcal{L}(M_{p}),\mathcal{L}(\tilde M_{p})) \ \lesssim \ n^{-\frac{1}{2}+\delta},\tag{\emph{i}}
\end{equation}
and there is a constant $C>0$ not depending on $n$, such that the event
\begin{equation}\label{eqnbpart}
d_{\textup{K}}(\mathcal{L}(\tilde M_{p}),\mathcal{L}(M_{p}^{\star}|X)) \ \leq \ C\,n^{-\frac{1}{2}+\delta} \tag{\emph{ii}}
\end{equation}
holds with probability at least $1-C/n$.
\end{proposition} 
The next two subsections will give the proofs of $(i)$ and $(ii)$ respectively.

\subsection{Proof of Proposition~\ref{PROP:BOOT}$(i)$}\label{sec:propiproof}

Recall $l_n = \big\lceil n^{\delta/\kappa}\wedge p\big\rceil$ in Theorem~\ref{THM:BOOT}, and define the integer 
$$k_n= l_n^2\wedge p.$$ 
We proceed by considering the bound
\[ \dK \big( \mathcal{L} (M_{p}) , \mathcal{L} (\tilde{M}_{p}) \big) \ \leq \ \I_n + \II_n + \III_n ,\]
where the terms on the right are defined by
\begin{align}
\I_n & \ = \ \dK \big( \mathcal{L} (M_{p}), \mathcal{L} (M_{k_n}) \big) \\[0.2cm]
\II_n & \ = \ \dK \big( \mathcal{L} (M_{k_n}), \mathcal{L} (\tilde{M}_{k_n}) \big)\label{eqn:IIdef} \\[0.2cm]
\III_n & \ = \ \dK \big( \mathcal{L} (\tilde{M}_{k_n}), \mathcal{L} (\tilde{M}_{p}) \big) .
\end{align}

In the remainder of this subsection, we will focus on establishing the bound
\begin{equation}\label{eqn:Ibound}
\I_n \lesssim n^{-\frac{1}{2}+\delta}.
\end{equation}
This will be sufficient to prove Proposition~\ref{PROP:BOOT}$(i)$, because the quantity $\III_n$ can be bounded using the same argument as for $\I_n$, and because the bound 
\begin{equation}\label{eqn:iibound}
\II_n\lesssim n^{-\frac{1}{2}+\delta}
\end{equation}
follows from Lemma~\ref{lem:II} of Section~\ref{sec:background}.

Going forward, including in Section~\ref{sec:ii}, we may assume without loss of generality that $k_n<p$, because if $p\leq k_n$, then $p=k_n$ and the quantities $\I_n$ and $\III_n$ become identically 0. Furthermore, the argument for bounding $\II_n$ does not depend on the relationship between $k_n$ and $p$. To proceed, first notice that the definition of $\I_n$ gives
\begin{equation*}
\begin{split}
\I_n & \ = \ \sup_{t \in \R} \, \Big| \P \Big( \max_{1 \leq j \leq p} S_{n,j}/\sigma_j^{\tau} \leq t \Big) - \P \Big( \max_{j \in J (k_n)} S_{n,j}/\sigma_j^{\tau} \leq t \Big) \Big|\\[0.2cm]
& \ = \ \sup_{t \in \R} \,  \P \big( A(t) \cap B (t) \big) ,
\end{split}
\end{equation*}
where we define the following events for arbitrary $t\in\R$,
\[ A(t) = \Big\{ \max_{j \in J (k_n)} S_{n,j}/\sigma_j^{\tau} \leq t \Big\} \quad \text{and} \quad B (t) = \Big\{ \max_{j \in J (k_n)^c} S_{n,j}/\sigma_j^{\tau} > t \Big\}. \]
For any real numbers $t_{1,n}$ and $t_{2,n}$ satisfying $t_{1,n} \leq t_{2,n}$, it can be checked that the following inclusion holds simultaneously for all $t\in\R$,
\[ A (t) \cap B(t) \ \subset \  A (t_{2,n}) \cup B (t_{1,n}) .\]
Therefore, a union bound gives
\begin{equation}\label{eqn:Inbound}
 \I_n 
\ \leq \ \P \big( A (t_{2,n})\big) +\P \big(B (t_{1,n}) \big).
\end{equation}

In order to choose suitable values of $t_{1,n}$ and $t_{2,n}$, define the parameter 
$$\omega=2\delta^2/\kappa,$$
which satisfies $\omega\in (0,\delta)$ for any choice of $\delta\in (0,1/2)$.  Also, define
 the integer
\begin{equation*}
 d_n = \big\lfloor (\ts\frac{\,\omega^2}{4}\,{\tt{r}}(R(l_n)))\vee 1\big\rfloor. 
\end{equation*}
Note that since ${\tt{r}}(R(l_n))\leq l_n$, it is clear that the inequalities $d_n\leq l_n\leq k_n\leq p$ hold for all $n$. Also, using the assumptions of Theorem~\ref{THM:BOOT}, we have the following lower bound on $d_n$,
\begin{equation}\label{eqn:mnlower}
d_n \ \gtrsim \ {\tt{r}}(R(l_n)) \ = \ \ts\frac{l_n^2}{\|R(l_n)\|_F^2} \ \gtrsim \ l_n^{\delta} \ \gtrsim \  n^{\delta^2/\kappa},
\end{equation}
which will be used later.

To finish the proof of~\eqref{eqn:Ibound}, it suffices to exhibit values $t_{1,n}$ and $t_{2,n}$  that satisfy the following three conditions:
\begin{align}
 t_{1,n} & \ \leq \ t_{2,n} \text{ \ \ \ holds for all large $n$}\label{eqn:tcond}\\[0.2cm]
 \P(B(t_{1,n})) & \ \lesssim \ n^{-1}\label{eqn:Bcond}\\[0.2cm]
 \P(A(t_{2,n})) &  \ \lesssim \ n^{-1/2+\delta}\label{eqn:Acond}.
\end{align}
In order to specify such values of $t_{1,n}$ and $t_{2,n}$, we will take them to be of the form
\begin{align*}
t_{1,n} & = c_1 \cdot k_n^{- \gamma(1-\tau) } \cdot \log (n) \\
t_{2,n} & = c_2 \cdot l_n^{- \gamma(1-\tau) } \cdot \sqrt{\log (d_n)},
\end{align*}
for certain constants $c_1,c_2>0$ that do not depend on $n$. First, the condition~\eqref{eqn:tcond} holds for any fixed choices of $c_1$ and $c_2$, due to the definitions of $k_n$, $l_n$, and $d_n$. 

Second, to establish \eqref{eqn:Bcond},  %
 define the parameter $q = \max \{ 2 / (\gamma(1-\tau)), \log (n), 3 \}$. For any $t > 0$, we have the tail bound
\begin{equation}\label{eqn:qnormbound}
 \P \big( B(t) \big) \leq t^{-q} \ \Big\| \max_{j \in J (k_n)^c} S_{n,j} / \sigma_j^\tau \Big\|_q^q ,
 \end{equation}
where $\|\cdot\|_q$ denotes the usual $L_q$ norm of a random variable. Due to Lemma~\ref{lem:qnorm} in Section~\ref{sec:background}, we have
$\| \frac{1}{\sigma_j} S_{n,j} \|_q \leq c q$ for all $j=1,\dots,p$, and so
\begin{align*}
\Big\| \max_{j \in \mathcal{J} (k_n)^c} S_{n,j} / \sigma_j^\tau \Big\|_q^q 
& \ \leq \ \sum_{j\in J(k_n)^c}\|S_{n,j}/\sigma_j^{\tau}\|_q^q\\[0.2cm]
& \ \leq \ (c q)^q \sum_{j \in J (k_n)^c} \sigma_j^{q(1-\tau)} \\[0.2cm]
& \ \lesssim \ (cq)^q \sum_{j = k_n+1}^p j^{-q \gamma(1-\tau)} \\[0.2cm]
& \ \lesssim \ \frac{(cq)^q}{q \gamma(1-\tau) - 1} \ k_n^{- q \gamma(1-\tau) + 1}.
\end{align*}
One can check that $q \asymp \log (n)$ and $\frac{c}{(q \gamma(1-\tau) - 1)^{1/q}} k_n^{1/q} \lesssim 1$. Therefore, if we take
$$t=e \cdot \Big\| \max_{j \in \mathcal{J} (k_n)^c} S_{n,j} / \sigma_j^\tau \Big\|_q$$
in~\eqref{eqn:qnormbound}, then there is a choice of $c_1$ for which $t_{1,n}$ satisfies
$$t\leq t_{1,n},$$
and furthermore,
$$\P(B(t_{1,n})) \ \leq \ \P(B(t)) \ \leq e^{-q} \ \leq \ n^{-1},$$
as needed for~\eqref{eqn:Bcond}.

Third, we now turn to~\eqref{eqn:Acond}. Observe that
\begin{equation}\label{eqn:beforeY}
\begin{split}
\P (A (t_{2,n})) & \ \leq \ \P \Big( \max_{j \in  J (k_n)} \tilde{S}_{n,j}/\sigma_j^{\tau} \, \leq  \, t_{2,n} \Big) + \II_n\\[0.2cm]
& \ \lesssim \ \P \Big( \max_{j \in J (l_n)} \tilde{S}_{n,j}/\sigma_j^{\tau} \, \leq \, t_{2,n} \Big) + n^{-1/2+\delta}
\end{split}
\end{equation}
where the second step follows from~\eqref{eqn:iibound} and the inclusion $J(l_n)\subset J(k_n)$. 

 As a preparatory step towards applying Theorem~\ref{thm:key} to the last line of~\eqref{eqn:beforeY}, we need the following basic observation. Let $(Y_j )_{j\in J(l_n)}$ be a generic set of random variables, and let $(a_j)_{j\in J(l_n)}$ be positive numbers satisfying the condition $\max_{j \in J(l_n)} a_j \leq b$, for some number $b$. It is straightforward to check that
\begin{equation}\label{eqn:tempybound}
 \P \Big( \max_{j \in J(l_n)} Y_j \, \leq \, t_{2,n} \Big) \ \leq \ \P \Big( \max_{j \in J (l_n)} a_j Y_j \, \leq \, b \, t_{2,n} \Big) .
 \end{equation}
 
Based on Assumption~\ref{assump1}, there is a positive constant $c_0$ not depending on $n$ such that the inequality $\sigma_j^{-(1-\tau)} \leq l_n^{\gamma(1-\tau)} / c_0$ holds for all $j \in J (l_n)$. Accordingly, we will use~\eqref{eqn:tempybound} with the choices $a_j = \sigma_j^{-(1-\tau)}$,  $b =l_n^{\gamma(1-\tau)} / c_0$, and $Y_j = \tilde{S}_{n,j}/\sigma_j^{\tau}$. Also, we may choose $c_2$ in the definition of $t_{2,n}$ to have the value $c_2=\omega c_0 \sqrt{2 (1- \omega)}$, which implies $b \, t_{2,n} = \omega\sqrt{2(1-\omega) \log (d_n)}$. Under these choices, the inequality~\eqref{eqn:tempybound} becomes %
\begin{equation*}
 \P \Big( \max_{j \in J (l_n)} \tilde{S}_{n,j}/\sigma_j^{\tau}\ \leq \ t_{2,n} \Big) \ \leq \ \P \bigg( \max_{j \in J (l_n)} \tilde{S}_{n,j} / \sigma_j \, \leq \, \omega\sqrt{2(1-\omega) \log (d_n)} \bigg).
 \end{equation*}
 \normalsize

Now, we apply Theorem~\ref{thm:key} to the right side, with $(l_n,d_n,\omega,\omega)$ playing the roles of $(N,k,\e,\delta)$ in that statement of that result. (Under these choices, the application of Theorem~\ref{thm:key} is justified because $\omega \in (0,1)$ and the inequalities $d_n\leq (\omega^2/4){\tt{r}}(R(l_n))$ and $d_n\geq 2$ hold for all large $n$.) Hence, 
$$  \P \Big( \max_{j \in J (l_n)} \tilde{S}_{n,j}/\sigma_j^{\tau}\ \leq \ t_{2,n} \Big)  \ \lesssim \ d_n^{\frac{-(1-\omega)(1-\omega)^2}{\omega}}(\log(d_n))^{\frac{1-\omega(2-\omega)-\omega}{2\omega}}.$$
\normalsize
Furthermore, using the lower bound on $d_n$ from~\eqref{eqn:mnlower}, there is a constant $c$ not depending on $n$ such that
\small
\begin{equation*}
\begin{split}
 \P \Big( \max_{j \in J (l_n)} \tilde{S}_{n,j}/\sigma_j^{\tau}\ \leq \ t_{2,n} \Big) 
&  \ \lesssim \ \big(n^{\delta^2/\kappa}\big)^{\frac{-(1-\omega)^3}{\omega}}\cdot (\log(n))^c\\[0.2cm]
& \ = \ n^{-\frac{(1-\omega)^3}{2}}\cdot (\log(n))^c\\[0.2cm]
& \ \lesssim \ n^{-\frac{1}{2}+3\omega}\\[0.2cm]
& \ \lesssim \ n^{-\frac{1}{2}+\delta}
\end{split}
\end{equation*}
\normalsize
which completes the proof.\qed

~\\

\subsection{Proof of Proposition~\ref{PROP:BOOT}$(ii)$}\label{sec:ii}

The proof of part $(ii)$ is structured mostly along the same lines as the proof of part $(i)$, and so we only sketch out the main steps for the sake of brevity. The current proof will also continue to use the same notation. Consider the inequality
\[ d_{\mathrm{K}} \big( \mathcal{L} (\tilde{M}_{p}), \mathcal{L} (M_{p}^\star | X) \big) \ \leq \ \I'_n \ + \  \II'_n (X)  \ + \  \III'_n (X) ,\]
where we define
\begin{align}
\I'_n & \ = \ d_{\mathrm{K}} \big( \L (\tilde{M}_{p}), \L (\tilde{M}_{k_n}) \big) \\[0.2cm]
\II'_n (X) & \ = \ d_{\mathrm{K}} \big( \L (\tilde{M}_{k_n}), \L (M_{k_n}^\star | X) \big)\label{eqn:IIprimedef} \\[0.2cm]
\III'_n (X) & \ = \ d_{\mathrm{K}} \big( \L (M_{k_n}^\star | X), \L (M_{p}^\star | X) \big) .
\end{align}
Note that $\I_n'$ is non-random, whereas $\II_n'(X)$ and $\III_n'(X)$ are random. 

Establishing a bound on $\I'_n$ of order $n^{-1/2+\delta}$ requires  no further work, because $\I_n'$ is equal to $\III_n$ in the proof of part $(i)$. Next, it follows from Lemma~\ref{lem:II} in Section~\ref{sec:background} that the event
\begin{equation}\label{eqn:IInprimebound}
\II_n'(X) \ \leq \ C\, n^{-\frac{1}{2}+\delta}
\end{equation}
holds with probability at least $1-C/n$. So, it remains to establish a bound on $\III_n'(X)$ of order $n^{-1/2+\delta}$.  For this purpose, let $t_{1,n}'$ and $ t_{2,n}'$ be any real numbers satisfying $t_{1,n}' \leq t_{2,n}'$. The reasoning leading up to~\eqref{eqn:Inbound} can be re-used to show that the following bound holds almost surely 
\begin{equation}\label{eqn:2starprobs}
 \mathrm{\III}_n' (X) \ \leq \ \P (A' (t_{2,n}') | X )  \ + \  \P (B' (t_{1,n}' ) | X) ,
 \end{equation}
where we define the following events for arbitrary $t\in\R$,
\[ A' (t) = \Big\{ \max_{j \in J (k_n)} S^\star_{n,j} \leq t \Big\} \quad \text{ \ \ and \ \ } \quad B' (t) = \Big\{ \max_{j \in J (k_n)^c} S^\star_{n,j} > t \Big\}. \]
To bound the two probabilities on the right side of~\eqref{eqn:2starprobs}, we will use values $t_{1,n}'$ and $t_{2,n}'$ having the form
\begin{align*}
t_{1,n}' & = c_1' \cdot k_n^{- \gamma(1-\tau) } \cdot (\log (n))^{3/2} \\[0.2cm]
t_{2,n}' & = c_2' \cdot l_n^{- \gamma(1-\tau)} \cdot \sqrt{\log (d_n)},
\end{align*}
for certain constants $c_1', c_2' > 0$ that do not depend on $n$.  (Note that for any such choices, the condition $t'_{1,n} \leq t'_{2,n}$ will hold for all large $n$.)

The probability $\P (A' (t_{2,n}') | X ) $ can be handled using the argument after~\eqref{eqn:beforeY} in the proof of part $(i)$, together with the bound~\eqref{eqn:IInprimebound}. Specifically, it can be shown that there are constants $C$ and $c_2'$ such that
\begin{align*}
\P (A' (t_{2,n}') | X ) 
& \ \leq \ \P \bigg( \max_{j \in J (l_n)} \tilde{S}_{n,j} \ \leq \ t_{2,n}' \bigg)  \ + \  \II_n' (X) \\[0.2cm]
& \ \leq \ C n^{-\frac{1}{2} + \delta} 
\end{align*}
holds with probability at least $1-C/n$. Finally, an argument analogous to the one used to establish~\eqref{eqn:Bcond} earlier shows there are constants $C$ and $c_1'$ such that 
$$\P (B' (t_{1,n}' ) | X) \leq n^{-1}$$ 
holds with probability at least $1-C/n$. (A more detailed version of this argument can be found in the proof of Lemma C.1 part (b) in~\citep{Lopes:2020}.) \qed

\subsection{Background results}\label{sec:background}
The following two bounds~\eqref{eqn:iiboundback} and~\eqref{eqn:iiprimeboundback} follow from the proofs of Propositions B.1 and C.1 in~\citep{Lopes:2020}. These proofs can be applied in the setting of Theorem~\eqref{THM:BOOT} in this paper with no essential changes.
\begin{lemma}\label{lem:II}
Fix any $\delta\in (0,1/2)$, and suppose the conditions of Theorem~\ref{THM:BOOT} hold. Also, let $\II_n$ and $\II_n'(X)$ be as defined in~\eqref{eqn:IIdef} and~\eqref{eqn:IIprimedef} respectively. Then,
\begin{equation}\label{eqn:iiboundback}
        \II_n \ \lesssim \ n^{\frac{1}{2}+\delta},
\end{equation}
and there is a constant $C>0$ not depending on $n$ such that the event
\begin{equation}\label{eqn:iiprimeboundback}
    \II_n'(X) \ \leq \ C n^{-\frac{1}{2}+\delta}
\end{equation}
holds with probability at least $1-C/n$.
\end{lemma}
 For the next lemma, recall that $\|\cdot\|_q$ denotes the $L_q$ norm of a random variable. This lemma is effectively a restatement of Lemma D.4 in~\citep{Lopes:2020}, and the proof given there can be applied in the same manner under the conditions used here.
\begin{lemma}\label{lem:qnorm}
Suppose the conditions of Theorem~\ref{THM:BOOT} hold, and define the parameter $q = \max \{ 2 / (\gamma(1-\tau)), \log (n), 3 \}$. Then, 
$$\max_{1\leq j\leq p}\| \ts\frac{1}{\sigma_j} S_{n,j} \|_q \ \lesssim \ q.$$
\end{lemma}

\references

\end{document}